\documentclass[a4paper]{amsart}
\usepackage{graphicx}
\usepackage{amssymb}
\usepackage{amsmath}
\usepackage{amsthm}
\usepackage{amscd}
\usepackage[all,2cell]{xy}

\UseAllTwocells \SilentMatrices
\newtheorem{thm}{Theorem}[section]

\newtheorem{lem}[thm]{Lemma}

\newtheorem{prop}[thm]{Proposition}
\theoremstyle{definition}

\theoremstyle{remark}
\newtheorem{rem}[thm]{\bf Remark}
\numberwithin{equation}{section}

\begin{document}
\title[A recollement of vector bundles]
{A recollement of vector bundles}
\author[Xiao-Wu Chen]{Xiao-Wu Chen}
\thanks{The author is supported by Special Foundation of President of The Chinese Academy of Sciences 
(No.1731112304061) and  National Natural Science Foundation of China (No.10971206)}
\date{Nov. 22, 2010}
\thanks{E-mail: xwchen$\symbol{64}$mail.ustc.edu.cn}

\maketitle

\dedicatory{}%
\commby{}%
\begin{center}
\end{center}

\begin{abstract}
For a weighted projective line, the stable category of its vector bundles modulo lines bundles
has a natural triangulated structure. We prove
that, for any positive integers $p, q, r$ and $r'$ with $r'\leq r$, there is an explicit
recollement of the stable category of vector bundles on a weighted projective line
of weight type $(p, q, r)$ relative to the ones on weighted projective lines of weight types $(p, q, r')$
and $(p, q, r-r'+1)$.
\end{abstract}

\section{Introduction}

Let $k$ be a field and  let $p, q, r$ be arbitrary positive integers. Denote by
$\mathbb{X}(p, q, r)$ the weighted projective line of weight type $(p, q,r)$ in the sense of
Geigle and Lenzing \cite{GL87}. We denote by $\mathrm{coh}\; \mathbb{X}(p, q, r)$ the category of
coherent sheaves on $\mathbb{X}(p, q, r)$ and by ${\rm vect}\; \mathbb{X}(p, q, r)$ the full subcategory consisting of vector
bundles. Following \cite{KLM1} a sequence $\eta\colon 0\rightarrow \mathcal{E}'\rightarrow \mathcal{E}
\rightarrow \mathcal{E}''\rightarrow 0 $ in ${\rm vect}\; \mathbb{X}(p, q, r)$  is \emph{distinguished exact}
provided that ${\rm Hom}(\mathcal{L}, \eta)$ are exact for all line bundles $\mathcal{L}$
on $\mathbb{X}(p, q, r)$. Observe that a distinguished exact sequence is exact. With the class
of distinguished exact sequences the category ${\rm vect}\; \mathbb{X}(p, q, r)$
is an exact category in the sense of Quillen \cite{Qui73}. Moreover, this exact category is \emph{Frobenius}, that is,
it has enough projective objects and enough injective objects such that the class of
projective objects coincides with the class of injective objects. In this setting, an object in
${\rm vect}\; \mathbb{X}(p, q, r)$ is projective if and only if it is a direct sum of line bundles.
Then by \cite[Chapter I, Theorem 2.8]{Ha88} the corresponding stable category $\underline{\rm vect}\; \mathbb{X}(p, q, r)$ of
${\rm vect}\; \mathbb{X}(p, q, r)$ modulo line bundles has a natural triangulated structure.

Recently, the stable category  $\underline{\rm vect}\; \mathbb{X}(p, q, r)$
of vector bundles receives a lot of attention. It is closely related to the category of
(graded) Cohen-Macaulay modules and then to the (graded) singularity category in the sense of Buchweitz \cite{Buc87} and Orlov \cite{Or04}.
Moreover, a version of Orlov's trichotomy theorem \cite{Or09} implies that this stable category  is also
related to the bounded derived category $\mathbf{D}^b({\rm coh}\; \mathbb{X}(p, q, r))$; see \cite[Section 7]{Lende}.
More recently, Kussin, Lenzing and Meltzer  \cite{KLM3}  prove that the stable category
 $\underline{\rm vect}\; \mathbb{X}(p, q, r)$ is triangle equivalent to the stable category of
the 2-flag category  of graded modules over $k[t]/(t^r)$  such that  the lengths of the two flags
are given by $p-1$ and $q-1$, respectively. Here, $t$ is an indeterminant with degree $1$. This result generalizes
their previous result in \cite{KLM2}, which gives a surprising link between weighted projective
lines and the (graded) submodule category of nilpotent operators; also see \cite{Ch11}.
 The latter category is studied intensively
by Ringel and Schmidmeier in a series of papers \cite{RS06, RS08, RS08'}.
Let us remark that the triangulated category $\underline{\rm vect}\; \mathbb{X}(p, q, r)$ has
nice homological properties such as having a tilting object and being fractionally Calabi-Yau.
For details, we refer to \cite{KLM1, KLM2, KLM3}.

The aim of this paper is to prove the following recollement \cite{BBD} consisting of the stable categories
of vector bundles on weighted projective lines.

 \vskip 5pt

\noindent {\bf Theorem.} \emph{
Let $p, q, r$ and $r'$ be positive integers such that $r'\leq r$. Then there exists a recollement of
triangulated categories
\[\xymatrixrowsep{3pc} \xymatrixcolsep{2pc}\xymatrix{
  \underline{\rm vect}\; \mathbb{X}(p, q, r')  \;\ar[rr]|-{} &&\;  \underline{\rm vect}\; \mathbb{X}(p, q, r) \; \ar[rr]|-{}
\ar@<1.2ex>[ll]^-{}\ar@<-1.2ex>[ll]_-{}&&
\; \underline{\rm vect}\; \mathbb{X}(p, q, r-r'+1). \ar@<1.2ex>[ll]^-{}\ar@<-1.2ex>[ll]_-{}
}\]
}

\vskip 5pt

This result is given in Theorem \ref{thm:maintheorem}, where the six functors in the recollement above
are given explicitly. A part of this recollement is obtained in \cite{CK} from the viewpoint of expansions of abelian categories. Let us
point out that the results in \cite{KLM3} suggest that  the recollement obtained here might relate to
the recollements constructed in \cite{CL}.

The paper is organized as follows. In Section 2, we collect some basic facts on adjoint
pairs and recollements. In Section 3, we recall some known results on the homogeneous coordinate
algebras of weighted projective lines. In particular, the relation among vector bundles,
graded Cohen-Macaulay modules, and graded singularity categories is recalled. We construct three exact
functors on the categories of graded modules over the homogeneous coordinate algebras in Section 4. These
functors induce the corresponding functors on the category of vector bundles. We state and prove our
main result in Section 5.

\section{Adjoint functors and recollements}

In this section we collect several well-known facts on adjoint functors
and recollements.

Let $F\colon \mathcal{A}\rightarrow \mathcal{B}$ and $G\colon \mathcal{B}
\rightarrow \mathcal{A}$  be two additive functors between additive categories. The pair $(F, G)$
is an \emph{adjoint  pair} provided that there is a functorial isomorphism of abelian groups
\begin{align}\label{equ:adjoint}
{\rm Hom}_\mathcal{B}(FA, B)\simeq {\rm Hom}_\mathcal{A}(A, GB).
\end{align}
This isomorphism induces the \emph{unit} $\eta\colon {\rm Id}_\mathcal{A}\rightarrow GF$ and
the \emph{counit }$\varepsilon\colon FG\rightarrow {\rm Id}_\mathcal{B}$ of the adjoint pair, both of
which are natural transformations. Recall that the functor $F$ is fully faithful if and
only if the unit $\eta$ is an isomorphism. We refer to \cite[Chapter IV]{MacL} for details.

Let $\mathcal{A'}$ be a Serre subcategory of an abelian category $\mathcal{A}$. Denote
by $\mathcal{A}/\mathcal{A}'$ the quotient abelian category in the sense of Gabriel \cite{Ga62}.
Consider an exact functor $F\colon \mathcal{A}\rightarrow \mathcal{B}$ between abelian categories, and
two Serre subcategories $\mathcal{A}'\subseteq \mathcal{A}$, $\mathcal{B}'\subseteq \mathcal{B}$
with $F\mathcal{A}'\subseteq \mathcal{B}'$. Then there is a uniquely induced exact functor $\bar{F}\colon
\mathcal{A}/\mathcal{A}'\rightarrow \mathcal{B}/\mathcal{B}'$.

\begin{lem}\label{lem:quotientAbel}
Let $F\colon \mathcal{A}\rightarrow
\mathcal{B}$ be an exact functor between abelian categories, which has an exact right
adjoint $G$. Assume that $\mathcal{A}'\subseteq \mathcal{A}$ and $\mathcal{B}'\subseteq \mathcal{B}$
are Serre subcategories such that $F\mathcal{A}'\subseteq \mathcal{B}'$ and $G\mathcal{B}'\subseteq \mathcal{A}'$.
Then the induced functor $\bar{F}\colon
\mathcal{A}/\mathcal{A}'\rightarrow \mathcal{B}/\mathcal{B}'$ is left adjoint to the induced functor
$\bar{G}\colon
\mathcal{B}/\mathcal{B}'\rightarrow \mathcal{A}/\mathcal{A}'$. Moreover,  if $F$ is fully faithful, then so is
 $\bar{F}$.
\end{lem}

\begin{proof}
Observe that the unit $\eta\colon {\rm Id}_\mathcal{A}\rightarrow GF$ (\emph{resp.},
the counit $\varepsilon\colon FG\rightarrow {\rm Id}_\mathcal{B}$) induces naturally a natural transformation
$\bar{\eta}\colon {\rm Id}_{\mathcal{A}/\mathcal{A}'}
\rightarrow \bar{G}\bar{F}$ (\emph{resp.}, $\bar{\varepsilon}\colon
\bar{F} \bar{G} \rightarrow {\rm Id}_{\mathcal{B}/\mathcal{B}'}$). Then we apply
\cite[Chapter IV, Section 1, Theorem 2(v)]{MacL} to deduce the adjoint pair  $(\bar{F}, \bar{G})$.
 Moreover, the corresponding unit (\emph{resp.}, counit) is $\bar{\eta}$ (\emph{resp.}, $\bar{\eta}$).

 If the functor $F$ is fully faithful, then  $\eta\colon {\rm Id}_\mathcal{A}\rightarrow GF$
 is an isomorphism. It follows that the natural transformation  $\bar{\eta}\colon {\rm Id}_{\mathcal{A}/\mathcal{A}'}
\rightarrow \bar{G}\bar{F}$ is also an isomorphism. This implies
that $\bar{F}$ is fully faithful; see \cite[Chapter IV, Section 3, Theorem 1]{MacL}.
\end{proof}

Replacing abelian categories by triangulated categories and Gabriel quotient by Verdier quotient \cite{Ver77} in Lemma \ref{lem:quotientAbel},
one obtains a triangulated analogue of Lemma \ref{lem:quotientAbel}; see \cite[Lemma 1.2]{Or04}.

Let $\mathcal{A}$ be an abelian category. Denote by $\mathbf{K}^b(\mathcal{A})$ and $\mathbf{D}^b(\mathcal{A})$ the
bounded homotopy category  and  the  bounded derived category of $\mathcal{A}$, respectively. Recall
that $\mathbf{D}^b(\mathcal{A})$  is the Verdier quotient triangulated category of
$\mathbf{K}^b(\mathcal{A})$ by the subcategory consisting of acyclic complexes.
We will always identify $\mathcal{A}$ as the full subcategory of $\mathbf{D}^b(\mathcal{A})$ formed by stalk complexes concentrated at degree zero. For details, we refer to \cite{Ver77, Ha88}

Let $\mathcal{B}$ be another abelian categories. Let $F\colon \mathcal{A}\rightarrow
\mathcal{B}$ be an exact functor. Then the functor $F$ extends naturally to a
triangle functor $\mathbf{D}^b(F)\colon \mathbf{D}^b(\mathcal{A})\rightarrow \mathbf{D}^b(\mathcal{B})$. The following fact could be proved  directly similar as Lemma \ref{lem:quotientAbel}; see \cite[Lemma 3.3.1(1)]{CK}.

\begin{lem}\label{lem:derived}
Let $F\colon \mathcal{A}\rightarrow
\mathcal{B}$ be an exact functor between abelian categories which has an exact right
adjoint $G$. Then the pair $(\mathbf{D}^b(F), \mathbf{D}^b(G))$ is adjoint. Moreover, if $F$
is fully faithful, then so is $\mathbf{D}^b(F)$.
\end{lem}

\begin{proof}
Observe that the isomorphism (\ref{equ:adjoint}) extends to the bounded homotopy categories. Then we
apply the triangulated analogue of Lemma \ref{lem:quotientAbel}.
\end{proof}

Let $\mathcal{A}$ be an additive category. Recall that a  sequence $X\stackrel{i} \rightarrow Y\stackrel{d} \rightarrow Z$
in  $\mathcal{A}$ is a \emph{kernel-cokernel sequence} if $i={\rm Ker}\;  d$ and $d={\rm Cok}\; d$.
By an \emph{exact category} in the sense of Quillen \cite{Qui73} we mean an additive category with a chosen class
of kernel-cokernel sequences which satisfies certain axioms. For an exact category $\mathcal{A}$
the sequence in the chosen class is called a \emph{conflation}. For example, an abelian category
is naturally an exact category with conflations induced by short exact sequences. More generally,
an extension-closed subcategory of an abelian category  is an exact category in the same manner.
An additive functor $F\colon \mathcal{A}\rightarrow \mathcal{B}$ between two exact categories
is \emph{exact} provided that it sends conflations to conflations. For details, we refer to
\cite[Appendix A]{Ke90}.

An exact category $\mathcal{A}$ is \emph{Frobenius} provided that it has enough projective and enough injective
objects, and that the class of projective objects coincides with the class of injective objects. The \emph{stable category}
$\underline{\mathcal{A}}$ of a Frobenius category $\mathcal{A}$ is defined as follows: the objects are
the same as in $\mathcal{A}$, while for two objects $X, Y$ the Hom set ${\rm Hom}_{\underline{\mathcal{A}}}(X, Y)$ is the quotient
of ${\rm Hom}_\mathcal{A}(X, Y)$ modulo the subgroup formed by those morphisms that factor though
a projective object; the composition of morphisms in $\underline{\mathcal{A}}$ is induced by the one of $\mathcal{A}$.
The stable category $\underline{\mathcal{A}}$ has a natural triangulated structure; see \cite[Chapter I, Section 2]{Ha88} and \cite[1.2]{Ke90}.

Let  $F\colon \mathcal{A}\rightarrow \mathcal{B}$ be an exact functor between two Frobenius
categories which sends projective objects to projective objects. Then there is a uniquely induced
functor $\underline{F}\colon \underline{\mathcal{A}}\rightarrow \underline{\mathcal{B}}$, which is triangle functor by
\cite[Chapter I, Lemma 2.8]{Ha88}.

We observe the following fact.

\begin{lem}\label{lem:stable}
Let $F\colon \mathcal{A}\rightarrow \mathcal{B}$ be an exact functor between two Frobenius categories which
sends projective objects to projective objects. Assume that $F$ admits a right adjoint
$G\colon \mathcal{B}\rightarrow \mathcal{A}$ which is also
exact. Then we have the following statements:
\begin{enumerate}
\item the functor $G$ sends projective objects to projective objects;
\item the pair $(\underline{F}, \underline{G})$ is adjoint;
\item if the functor $F$ is fully faithful, so is $\underline{F}$.
\end{enumerate}
\end{lem}

\begin{proof}
(1) follows from a general fact that a right adjoint of an exact functor preserves injective
objects. (2) and (3) are easy, which could be proved by the same argument as in the proof
of Lemma \ref{lem:quotientAbel}.
\end{proof}

Recall that a diagram of triangle functors between triangulated
categories
\begin{equation*}\label{eq:recollement}
\xymatrixrowsep{3pc} \xymatrixcolsep{2pc}\xymatrix{
  \mathcal{T}'\;\ar[rr]|-{i}&&\;\mathcal{T}\; \ar[rr]|-{j}
  \ar@<1.2ex>[ll]^-{i_\rho}\ar@<-1.2ex>[ll]_-{i_\lambda}&&
  \;\mathcal{T}''\ar@<1.2ex>[ll]^-{j_\rho}\ar@<-1.2ex>[ll]_-{j_\lambda}
}\end{equation*} forms a \emph{recollement} \cite{BBD}, provided that the
following conditions are satisfied:
\begin{enumerate}
\item[{\rm (R1)}] The pairs $(i_\lambda, i)$, $(i, i_\rho)$, $(j_\lambda, j)$, and $(j,j_\rho)$ are adjoint.
\item[{\rm (R2)}] The functors $i$, $j_\lambda$, and $j_\rho$ are fully faithful.
\item[{\rm (R3)}] ${\rm Im}\; i={\rm Ker}\;  j$.
\end{enumerate}
Here for an additive functor $F$, ${\rm Im}\; F$ and ${\rm Ker}\; F$ denotes the essential
image and kernel of $F$, respectively. Recall that in this situation $j$ induces a triangle equivalence
$\mathcal{T}/{{\rm Ker}\; j}\simeq \mathcal{T}''$, where $\mathcal{T}/{{\rm Ker}\; j}$ denotes
the Verdier quotient category \cite{Ver77}.

The following two results are well known.

\begin{lem}\label{lem:recollement1}
Let $i\colon \mathcal{T}'\rightarrow \mathcal{T}$ be a fully faithful triangle functor which
admits a left adjoint $i_\rho$ and a right adjoint $i_\rho$. Then we have a recollement
of triangulated categories
\begin{equation*}\label{eq:recollement}
\xymatrixrowsep{3pc} \xymatrixcolsep{2pc}\xymatrix{
  \mathcal{T}'\;\ar[rr]|-{i}&&\;\mathcal{T}\; \ar[rr]|-{q}
  \ar@<1.2ex>[ll]^-{i_\rho}\ar@<-1.2ex>[ll]_-{i_\lambda}&&
  \;\mathcal{T}/{{\rm Im}\; i}\ar@<1.2ex>[ll]^-{}\ar@<-1.2ex>[ll]_-{}
}\end{equation*}
where $q\colon \mathcal{T}\rightarrow \mathcal{T}/{{\rm Im}\; i}$ denotes the quotient functor.
\end{lem}

\begin{proof}
Observe that the functors $i_\rho$ and $i_\lambda$ are triangle functors; see \cite[Lemma 8.3]{Ke96}.
The remaining follows directly from \cite[Propositions 1.5 and 1.6]{BK}.
\end{proof}

Recall that a \emph{thick} subcategory of a triangulated category $\mathcal{T}$ means
a full triangulated subcategory which is closed under taking direct summands. For a class $\mathcal{S}$
of objects in $\mathcal{T}$, denote by ${\rm thick} \langle \mathcal{S} \rangle$ the smallest thick subcategory
of $\mathcal{T}$ containing $\mathcal{S}$, which is called the thick subcategory \emph{generated} by $\mathcal{S}$;
compare \cite[p.70]{Ha88}.

\begin{lem}\label{lem:recollement2}
Suppose that we are given a diagram of triangle functors satisfying (R1) and (R2).
\begin{equation*}\label{eq:recollement}
\xymatrixrowsep{3pc} \xymatrixcolsep{2pc}\xymatrix{
  \mathcal{T}'\;\ar[rr]|-{i}&&\;\mathcal{T}\; \ar[rr]|-{j}
  \ar@<1.2ex>[ll]^-{i_\rho}\ar@<-1.2ex>[ll]_-{i_\lambda}&&
  \;\mathcal{T}''\ar@<1.2ex>[ll]^-{j_\rho}\ar@<-1.2ex>[ll]_-{j_\lambda}
}\end{equation*}
Assume that $j\circ i\simeq 0$ and ${\rm thick}\langle {\rm Im}\; i\cup {\rm Im}j_\lambda\rangle=\mathcal{T}$. Then this diagram
of functors is a recollement.
\end{lem}

\begin{proof}
It suffices to show that an object $X$ in $\mathcal{T}$ satisfying $jX\simeq 0$ lies
in ${\rm Im}\; i$. Consider the triangle $X'\rightarrow X\rightarrow ii_\lambda X\rightarrow X'[1]$,
where $X\rightarrow ii_\lambda X$ is given by the unit of the adjoint pair  $(i_\lambda, i)$ and $[1]$ denotes
the translation functor of $\mathcal{T}$. It follows that $jX'\simeq 0$ and ${\rm Im}\; i\subseteq X'^\perp$; here,
$X'^\perp=\{Y\in \mathcal{T}\; |\;{\rm Hom}_\mathcal{T}(X', Y[n])=0, \; n\in \mathbb{Z}\}$ which
is a thick subcategory. Here, for each $n\geq 1$, [n] (\emph{resp.}, [-n]) denotes the
 $n$-th power of the translation functor [1] (\emph{resp.}, the inverse [-1] of the translation functor). Observe that ${\rm Im}\; j_\lambda\subseteq X'^\perp$ by the adjoint pair $(j_\lambda, j)$.
Then it follows from  ${\rm thick}\langle {\rm Im}\; i\cup {\rm Im}j_\lambda\rangle=\mathcal{T}$ that $\mathcal{T}\subseteq X'^\perp$, which forces that $X'\simeq 0$. Then we
have $X\simeq ii_\lambda X$. We are done.
\end{proof}

\section{Homogenous coordinate algebras}

In this section we recall some basic facts on the homogeneous coordinate algebras of
weighted projective lines. In particular, the relation among vector bundles on weighted
projective lines, graded Cohen-Macaulay modules and the graded singularity category of
the homogeneous coordinate algebras is recalled.

Let $\mathbf{p}=(p_1, p_2, \cdots, p_n)$ be a sequence of positive integers with $n\geq 2$, which
is called a \emph{weight sequence}. Denote by $\mathbf{L}=\mathbf{L}(\mathbf{p})$
the rank one abelian group generated by $\vec{x}_1, \vec{x}_2, \cdots, \vec{x}_n$
subject to the relations $p_1\vec{x}_1=p_2\vec{x}_2= \cdots =p_n\vec{x}_n$.
The torsionfree element $\vec{c}=p_1\vec{x}_1$ in $\mathbf{L}$ is called the \emph{canonical element}.
Recall that each element $\vec{l}$ in $\mathbf{L}$ can be uniquely expressed in this
\emph{normal form} $\vec{l}=\sum_{i=1}^n l_i \vec{x}_i+l \vec{c}$ such that
$l\in \mathbb{Z}$ and $0\leq l_i< p_i$ for each $i$; see \cite[1.2]{GL87}. In what follows,
 all elements in $\mathbf{L}$ will be written in their normal forms.

Let $k$ be an arbitrary field. Denote by $\mathbb{P}_k^1$ the projective line over $k$.
For each rational point $\lambda$ of $\mathbb{P}_k^1$ we fix a choice of its homogeneous
coordinates $\lambda=[\lambda_0: \lambda_1]$. Let $\boldsymbol\lambda=(\lambda_1, \lambda_2,\cdots, \lambda_n)$
be a sequence of pairwise distinct rational points of $\mathbb{P}_k^1$, which is called a \emph{parameter sequence}.

Denote by $\mathbb{X}(\mathbf{p}, \boldsymbol\lambda)$ the \emph{weighted projective line} \cite{GL87} with weight
sequence $\mathbf{p}$ and parameter sequence $\boldsymbol\lambda$. Recall
that $\mathbf{S}=\mathbf{S}(\mathbf{p}, \boldsymbol\lambda)$ the \emph{homogeneous coordinate algebra} of
$\mathbb{X}(\mathbf{p}, \boldsymbol\lambda)$ is defined by .
$$\mathbf{S}(\mathbf{p}, \boldsymbol\lambda)=k[U, V, X_1, X_2, \cdots, X_n]/(X_i^{p_i}+\lambda_{i1}U-\lambda_{i0}V,\; 1\leq i\leq n).$$
Here, we recall that $\lambda_i=[\lambda_{i0}:\lambda_{i1}]$. We write $u$, $v$
and  $x_i$ for the canonical image of  $U$, $V$ and $X_i$ in $\mathbf{S}$, $1\leq i\leq n$.
The algebra $\mathbf{S}$ is naturally $\mathbf{L}$-graded by means of
${\rm deg}\; u={\rm deg}\; v=\vec{c}$ and ${\rm deg}\; x_i = \vec{x}_i$.
Observe that $\mathbf{S}$ is (graded) noetherian.

We denote by ${\rm mod}^\mathbf{L}\; \mathbf{S}$ the abelian category of finitely generated
$\mathbf{L}$-graded $\mathbf{S}$-modules. A graded $\mathbf{S}$-module is written as
$M=\oplus_{\vec{l}\in \mathbf{L}} M_{\vec{l}}$, where $M_{\vec{l}}$ is the homogeneous
component of degree $\vec{l}$. For an element $\vec{l}$ in $\mathbf{L}$ the
\emph{shifted module} $M(\vec{l})$ is the same as $M$ as ungraded $S$-modules, while it
is graded such that $M(\vec{l})_{\vec{l}'}=M_{\vec{l}+\vec{l}'}$. This yields the \emph{degree-shift}
functor $(\vec{l})\colon {\rm mod}^\mathbf{L}\; \mathbf{S}\rightarrow {\rm mod}^\mathbf{L}\; \mathbf{S}$, which
is clearly an automorphism of categories. Observe that a complete set of representatives
of pairwise non-isomorphic indecomposable projective modules in ${\rm mod}^\mathbf{L}\; \mathbf{S}$ is given by
$\{\mathbf{S}(\vec{l})\; |\; \vec{l}\in \mathbf{L}\}$.
Here, we view $\mathbf{S}$ as a graded $\mathbf{S}$-module generated by its homogeneous component
of degree zero.

We remark that the $\mathbf{L}$-graded algebra $\mathbf{S}(\mathbf{p}, \boldsymbol \lambda)$, even up to isomorphism, might depend on the choice
of the homogeneous coordinates of the parameters $\lambda_i$'s.  However, the category ${\rm mod}^\mathbf{L}\; (\mathbf{p}, \boldsymbol \lambda)$  of finitely generated $\mathbf{L}$-graded $\mathbf{S}(\mathbf{p}, \boldsymbol \lambda)$-modules, up to equivalence,  does not
depend on such choice.

We observe that  the algebra $\mathbf{S}=\mathbf{S}(\mathbf{p}, \boldsymbol \lambda)$ is \emph{graded local}, that is, it has a unique maximal
homogeneous ideal $\mathfrak{m}=(x_1, x_2, \cdots, x_n)$. Consider $k=\mathbf{S}/\mathfrak{m}$ the \emph{trivial
module} of $\mathbf{S}$, which is concentrated at degree zero.
Then the set $\{k(\vec{l})\; |\; \vec{l}\in \mathbf{L}\}$ is a complete
set of representatives of pairwise non-isomorphic graded simple $\mathbf{S}$-modules.

\begin{lem}\label{lem:basicofS}
Use the notation above. Then the following statements hold:
\begin{enumerate}
\item the graded $\mathbf{S}$-module $\mathbf{S}$ has injective dimension two, in particular, the
algebra $\mathbf{S}$ is graded Gorenstein;
\item the algebra $\mathbf{S}$ is an graded isolated singularity, that is, for each
homogeneous non-maximal prime ideals $\mathfrak{p}$ the homogeneous localization $\mathbf{S}_\mathfrak{p}$
has finite graded global dimension.
\end{enumerate}
\end{lem}

\begin{proof}
(1) follows from the observation that $\mathbf{S}$ has (graded) Krull dimension two and
it is a complete intersection, and then Gorenstein; compare \cite[Proposition 1.3]{GL87}.
(2) follows from \cite[1.6]{GL87}.
\end{proof}

Denote by $k[u, v]$ the polynomial ring with two variables which is $\mathbf{L}$-graded such
that ${\rm deg}\; u={\rm deg}\; v=\vec{c}$. The following embedding of $\mathbf{L}$-graded algebras is known as
the \emph{core homomorphism}
$$k[u, v]\longrightarrow \mathbf{S},$$
which sends $u$ and $v$ to $u$ and $v$, respectively. Observe that the algebra $\mathbf{S}$ is a
finitely generated free module over $k[u, v]$ via this core homomorphism.
More precisely, each homogeneous component $\mathbf{S}_{\vec{l}}$ has an explicit basis
$\{ \prod_{i=1}^n x_i^{l_i} u^av^b\; | \; a+b=l, a, b\geq 0 \}$.
Hence the $k[u, v]$-module $\mathbf{S}$ has a  homogeneous basis
$\{ \prod_{i=1}^n x_i^{l_i} \; | \; 0\leq l_i < p_i,\; 1\leq i\leq n\}$.
For details, we refer to  the proof of \cite[Proposition 1.3]{GL87}.
Observe that again via the core homomorphism, a graded $\mathbf{S}$-module $M$ induces graded $k[u, v]$-modules
$M|_{\vec{l}+\mathbb{Z}\vec{c}}$ for all $\vec{l}\in \mathbf{L}$.  Here $M|_{\vec{l}+\mathbb{Z}\vec{c}}=\oplus_{\vec{l}'\in \vec{l}+\mathbb{Z}\vec{c}} M_{\vec{l}'}$.

The following lemma is easy.

\begin{lem}\label{lem:f.g.}
A graded $\mathbf{S}$-module $M$ is finitely generated if and only if all
the induced $k[u, v]$-modules $M|_{\vec{l}+\mathbb{Z}\vec{c}}$ are finitely
generated. \hfill $\square$
\end{lem}

Recall that a module $M$ in ${\rm mod}^\mathbf{L}\; \mathbf{S}$ is called
\emph{(maximal) Cohen-Macaulay} provided that
${\rm Ext}^i_{{\rm mod}^\mathbf{L}\; \mathbf{S}}(M, \mathbf{S}(\vec{l}))=0$ for all
$i\geq 1$ and $\vec{l}$ in $\mathbf{L}$. Denote by ${\rm CM}^\mathbf{L}(\mathbf{S})$
the full subcategory consisting of Cohen-Macaulay modules. Observe that projective
$\mathbf{S}$-modules are Cohen-Macaulay and that ${\rm CM}^\mathbf{L}(\mathbf{S})$ is
extension-closed in ${\rm mod}^\mathbf{L}\; \mathbf{S}$. Hence ${\rm CM}^\mathbf{L}(\mathbf{S})$
becomes naturally an exact category.  Since $\mathbf{S}$ is Gorenstein, this exact category
is Frobenius; moreover, an object $M$ in ${\rm CM}^\mathbf{L}(\mathbf{S})$ is projective
if and only if it is a projective $\mathbf{S}$-module; see \cite[Lemma 4.2.2]{Buc87}. We denote by $\underline{{\rm CM}}^\mathbf{L}(\mathbf{S})$ the stable category; it is a triangulated category.

\begin{lem}\label{lem:MCM}
A graded $\mathbf{S}$-module $M$ is Cohen-Macaulay if and only if all
the induced $k[u, v]$-modules $M|_{\vec{l}+\mathbb{Z}\vec{c}}$ are (finitely generated)
projective.
\end{lem}

\begin{proof}
Recall that the algebra $\mathbf{S}$ is graded Gorenstein of injective dimension two.
 Then by a graded version of local duality the module $M$ is Cohen-Macaulay if
and only if ${\rm Hom}_{{\rm mod}^\mathbf{L}\; \mathbf{S}}(k(\vec{l}), M)=0={\rm Ext}^1_{{\rm mod}^\mathbf{L}\; \mathbf{S}}(k(\vec{l}), M)$ for all $\vec{l}$ in $\mathbf{L}$. Then this is equivalent to that
${\rm Hom}_{{\rm mod}^\mathbf{L}\; k[u, v]}(k(\vec{l}), M)=0={\rm Ext}^1_{{\rm mod}^\mathbf{L}\; k[u, v]}(k(\vec{l}), M)$
for all $\vec{l}$ in $\mathbf{L}$; see the third paragraph of the proof of \cite[Theorem 5.1]{GL87}. Observe
that the algebra $k[u, v]$ has graded global dimension two. Then this is equivalent to that $M$ is
a graded projective $k[u, v]$-module.
\end{proof}

Denote by $\mathbf{D}^b({\rm mod}^\mathbf{L}\; \mathbf{S})$ the bounded derived category of
${\rm mod}^\mathbf{L}\; \mathbf{S}$. We identify ${\rm mod}^\mathbf{L}\; \mathbf{S}$ as the full subcategory
of $\mathbf{D}^b({\rm mod}^\mathbf{L}\; \mathbf{S})$ consisting of stalk complexes concentrated
at degree zero. Denote by ${\rm perf}^\mathbf{L}(\mathbf{S})$ the full triangulated subcategory of
$\mathbf{D}^b({\rm mod}^\mathbf{L}\; \mathbf{S})$ consisting of \emph{perfect complexes}. Here we recall
that perfect complexes in $\mathbf{D}^b({\rm mod}^\mathbf{L}\; \mathbf{S})$ are those complexes isomorphic
to a bounded complex of finitely generated projective modules in ${\rm mod}^\mathbf{L}\; \mathbf{S}$. Following
\cite{Buc87} and \cite{Or04} the \emph{graded singularity category} of $\mathbf{S}$ is defined
to be  the Verdier quotient category ${\rm D}_{\rm sg}^\mathbf{L}(\mathbf{S})= \mathbf{D}^b({\rm mod}^\mathbf{L}\; \mathbf{S})/{{\rm perf}^\mathbf{L} (\mathbf{S}})$.

Consider the composite ${\rm mod}^\mathbf{L}\; \mathbf{S}\hookrightarrow \mathbf{D}^b({\rm mod}^\mathbf{L}\; \mathbf{S})
\stackrel{q}\rightarrow \mathbf{D}^\mathbf{L}_{\rm sg}(\mathbf{S})$, where the first functor identifies
a module $M$ as a stalk complex concentrated at degree zero, and $q$ denotes
the quotient functor. In this way $qM$ becomes an object in $\mathbf{D}^\mathbf{L}_{\rm sg}(\mathbf{S})$.
We restrict this composite functor to ${\rm CM}^\mathbf{L}(\mathbf{S})$. Observe
that $qP$ is isomorphic to zero for a projective $\mathbf{S}$-module $P$. Then we have an induced
functor $\underline{{\rm CM}}^\mathbf{L}(\mathbf{S})\rightarrow \mathbf{D}^\mathbf{L}_{\rm sg}(\mathbf{S})$.

\begin{lem}\label{lem:singularity}
Keep the notation above. Then  the following statements hold:
\begin{enumerate}
\item the induced functor $\underline{{\rm CM}}^\mathbf{L}(\mathbf{S})\rightarrow \mathbf{D}^\mathbf{L}_{\rm sg}(\mathbf{S})$
is a triangle equivalence;
\item  $\mathbf{D}^\mathbf{L}_{\rm sg}(\mathbf{S})={\rm thick}\langle qk(\vec{l})\; |\; \vec{l}\in \mathbf{L}\rangle $.
\end{enumerate}
\end{lem}

\begin{proof}
(1) Recall that the algebra $\mathbf{S}$ is graded Gorenstein. Then the triangle equivalence
follows from an $\mathbf{L}$-graded version of Buchweitz's theorem \cite[Theorem 4.4.1]{Buc87}.

(2) We use the fact that $\mathbf{S}$ is a graded isolated singularity; see Lemma \ref{lem:basicofS}(2). Then the statement
follows from an $\mathbf{L}$-graded version of \cite[Proposition A.2]{KMV08}; also see \cite[Corollary 2.4]{Ch11'} and
 compare \cite[Proposition 2.7]{Or10}.
\end{proof}

We denote by ${\rm coh}\; \mathbb{X}$ the abelian category of coherent sheaves on $\mathbb{X}=\mathbb{X}(\mathbf{p}, \boldsymbol \lambda)$. Recall that ${\rm coh}\; \mathbb{X}$ is equivalent to the quotient category ${\rm mod}^\mathbf{L}\; \mathbf{S}/ {{\rm mod}_0^\mathbf{L}\; \mathbf{S}}$, where ${\rm mod}^\mathbf{L}_0\; \mathbf{S}$ denotes the Serre subcategory consisting of
finite dimension modules. This corresponds to a quotient functor
${\rm mod}^\mathbf{L}\; S\rightarrow {\rm coh}\; \mathbb{X}$, which is known as the \emph{sheafification
functor} \cite{GL87}.  For details, see \cite[1.8]{GL87}.  The degree-shift functors $(\vec{l})$ yield the \emph{twist}
functors on ${\rm coh}\; \mathbb{X}$, which are still denoted by $(\vec{l})$.
Then the sheafification functor sends $\mathbf{S}(\vec{l})$ to the \emph{twisted structure
sheaf} $\mathcal{O}_\mathbb{X}(\vec{l})$ of $\mathbb{X}$.

 Recall that a locally free sheaf on $\mathbb{X}$ is called \emph{vector bundle}.
 A \emph{line bundle} is a vector bundle of rank one. Recall that a complete set of
representatives of pairwise non-isomorphic line bundles on $\mathbb{X}$ is given by $\{\mathcal{O}_\mathbb{X}(\vec{l})\; |\; \vec{l}\in \mathbf{L}\}$; see \cite[Proposition 2.1]{GL87}.

 We denote by ${\rm vect}\; \mathbb{X}$ the full subcategory of ${\rm coh}\; \mathbb{X}$ consisting of
vector bundles. This is an extension-closed subcategory and thus becomes an exact category. However,
we are not interested in this exact category, since it is not Frobenius in general.

Following \cite{KLM1}  a  sequence $\eta\colon 0\rightarrow \mathcal{E}'\rightarrow \mathcal{E}
\rightarrow \mathcal{E}''\rightarrow 0 $ in ${\rm vect}\; \mathbb{X}$  is \emph{distinguished exact}
provided that ${\rm Hom}(\mathcal{O}_\mathbb{X}(\vec{l}), \eta)$ are exact for all $\vec{l}$ in $\mathbf{L}$.
Observe that a distinguished exact sequence is exact in ${\rm coh}\; \mathbb{X}$. With the class
of distinguished exact sequences the category ${\rm vect}\; \mathbb{X}$ is a Frobenius category such that
an object  is projective if and only if it is a direct sum of line bundles. Denote by $\underline{\rm vect}\; \mathbb{X}$
the corresponding stable category, which is triangulated. This category is called the \emph{stable category
of vector bundles} on $\mathbb{X}$ \cite{KLM1}.

We have the following result in \cite{KLM1}. Recall that an exact functor between two exact category is
an \emph{equivalence of exact categories} provided that it is an equivalence of categories and its quasi-inverse
is also an exact functor.

\begin{lem}\label{lem:MCM-vb}
The sheafification functor induces an equivalence of exact categories ${\rm CM}^\mathbf{L}(S)\simeq {\rm vect}\; \mathbb{X}$,
which further induces a triangle equivalence $\underline{\rm CM}^\mathbf{L}(S)\simeq \underline{\rm vect}\; \mathbb{X}$.
\end{lem}

\begin{proof}
This follows from \cite[Theorem 5.1]{GL87} immediately.
\end{proof}

We observe that the degree-shift functors $(\vec{l})$ act on $\underline{\rm CM}^\mathbf{L}(\mathbf{S})$
and $\mathbf{D}_{\rm sg}^\mathbf{L}(\mathbf{S})$ naturally. Similarly, the twist functors
$(\vec{l})$ act on $\underline{\rm vect}\; \mathbb{X}$.

\begin{prop}\label{prop:vb-singularity}
Keep the notation as above. Then there is a triangle equivalence $\underline{\rm vect}\; \mathbb{X}\simeq
\mathbf{D}_{\rm sg}^\mathbf{L}(\mathbf{S})$, which is compatible with the degree-shift functors
and the twist functors.
\end{prop}

\begin{proof}
Combine Lemmas \ref{lem:singularity} and \ref{lem:MCM-vb}. Observe that the equivalences in the
two lemmas are compatible with the degree-shift functors and the twist functors.
\end{proof}

\section{Functors on graded modules and sheaves}

In this section we construct three functors on the graded module categories of the
homogeneous coordinate rings of weighted projective lines. These functors will induce
a recollement of the stable category of vector bundles. Let us point out that the construction
of these functors is essentially contained in \cite[Section 9]{GL90}; also see \cite[Section 4]{CK}.

Let $\mathbf{p}=(p_1, p_2, \cdots, p_n)$ the weight sequence in Section 3.
Fix a positive integer $p_n'\leq p_n$. Write $\mathbf{p}'=(p_1, p_2, \cdots, p'_n)$.
Denote $\mathbf{L}'=\mathbf{L}(\mathbf{p}')$. Consider the following
injective map $\phi'\colon \mathbf{L}'\rightarrow \mathbf{L}$ which sends an element $\vec{l}=\sum_{i=1}^n l_i
 \vec{x}_i+l\vec{c}$ to $\phi'(\vec{l})=\sum_{i=1}^n l_i
 \vec{x}_i+l\vec{c}$. Here, the element $\vec{l}$ in $\mathbf{L}'$ is in its normal form, that is,
$0\leq l_i<p'_i$ and $l\in \mathbb{Z}$, where $p'_i=p_i$ for $i<n$. Observe that in general, the map
$\phi'$ is not a homomorphism of groups; moreover, an element $\vec{l}$ in $\mathbf{L}$ lies in
the image of $\phi'$ if and only if $l_n<p'_n$.

Let $\boldsymbol \lambda$ be the parameter sequence in Section 3.
We denote by  $\mathbf{S}'=\mathbf{S}(\mathbf{p}', \boldsymbol \lambda)$ the homogeneous coordinate
 algebra of the weighted projective line $\mathbb{X}'=\mathbb{X}(\mathbf{p}', \boldsymbol \lambda)$.
 Then $\mathbf{S}'$ is naturally $\mathbf{L}'$-graded.

 We will define a functor $i'\colon {\rm mod}^{\mathbf{L}'}\; \mathbf{S}'\rightarrow {\rm mod}^\mathbf{L}\; \mathbf{S}$ as follows. For an  $\mathbf{L}'$-graded $\mathbf{S}'$-module $M$, define $i'M=\oplus_{\vec{l}\in \mathbf{L}} (i'M)_{\vec{l}}$ such that
$(i'M)_{\vec{l}}=M_{\phi'^{-1}(\vec{l}-l_n\vec{x}_n)}$ if $0\leq l_n< p_n-p'_n$, and $(i'M)_{\vec{l}}=M_{\phi'^{-1}(\vec{l}-(p_n-p'_n)\vec{x}_n)}$, otherwise. The action of $u$, $v$ and $x_i$ on
$i'M$ is induced by the one on $M$, except that $x_n$ acts as the identity on $(i'M)_{\vec{l}}$ provided
that $l_n< p_n-p'_n$. The obtained $\mathbf{L}$-graded $\mathbf{S}$-module $i'M$ is finitely generated by Lemma \ref{lem:f.g.}.
The action of $i'$ on morphisms is defined naturally.

\begin{lem}\label{lem:i'}
Use the notation above. Then the following statements hold:
\begin{enumerate}
\item the functor $i'$ is exact and fully faithful;
\item $i'(M(\vec{x}_i))=(i' M)(\vec{x}_i)$ for any graded $\mathbf{S}'$-module $M$ and $1\leq i< n$;
\item $i' (\mathbf{S}'(\vec{l}))\simeq \mathbf{S}(\phi'(\vec{l}))$ for all $\vec{l}\in \mathbf{L}'$;
\item $i' (k(\vec{l}))=k(\phi'(\vec{l}))$ if $l_n>0$, and $i'(k(\vec{l}))=\mathbf{S}(\phi'(\vec{l}))/(x_1, x_2, \cdots, x_{n-1}, x_{n}^{p_n-p'_n+1})$, otherwise;
\item if $M$ is a Cohen-Macaulay $\mathbf{S}'$-module, then  $i'M$ is a Cohen-Macaulay $\mathbf{S}$-module.
\end{enumerate}
\end{lem}

\begin{proof} The statements (1), (2) and (4) are direct from the construction of $i'$.
 The isomorphism in (3) is obtained by comparing the explicit bases of homogeneous components
 of the two graded $\mathbf{S}$-modules. The last statement follows from Lemma \ref{lem:MCM}.
\end{proof}

 We will define another functor $i'_\lambda \colon {\rm mod}^{\mathbf{L}}\; \mathbf{S}\rightarrow {\rm mod}^{\mathbf{L}'}\; \mathbf{S}'$ as follows. For an $\mathbf{L}$-graded $\mathbf{S}$-module $N$, we set $i'_\lambda N=\oplus_{\vec{l}\in \mathbf{L}'} (i'_\lambda N)_{\vec{l}}$ such that  $(i'_\lambda N)_{\vec{l}}=N_{\phi'(\vec{l})+(p_n-p'_n)\vec{x}_n}$ for all $\vec{l}$ in $\mathbf{L}'$. The action of $u$, $v$ and  $x_i$ on $i'_\lambda N$ is induced by the one on $N$, except that $x_n$ acts by
$x_n^{p_n-p'_n+1}$ (in $\mathbf{S}$) on $(i'_\lambda N)_{\vec{l}}$  provided that $l_n=p'_n-1$.
 The obtained $\mathbf{L}'$-graded $\mathbf{S}'$-module $i'_\lambda N$ is finitely generated by Lemma \ref{lem:f.g.}.
The action of $i'_\lambda$ on morphisms is defined naturally.

\begin{lem}\label{lem:leftadj}
Use the notation above. Then the following statements hold:
\begin{enumerate}
\item the functor $i'_\lambda$ is exact;
\item $i'_\lambda (N(\vec{x}_i))=(i'_\lambda N)(\vec{x}_i)$ for any graded $\mathbf{S}$-module $N$ and $1\leq i< n$;
\item $i'_\lambda (\mathbf{S}(\vec{l}))\simeq \mathbf{S}'(\phi'^{-1}(\vec{l}))$  if $0\leq l_n< p'_n$, and $i'_\lambda (\mathbf{S}(\vec{l}))\simeq \mathbf{S}'(\phi'^{-1}(\vec{l}-l_n\vec{x}_n)+\vec{c})$, otherwise;
\item $i'_\lambda (k(\vec{l}))=k(\phi'^{-1}(\vec{l}-\vec{x}_n)+\vec{x}_n)$ if $1\leq l_n\leq p'_n$, and $i'_\lambda (k(\vec{l}))=0$, otherwise;
\item if $N$ is a Cohen-Macaulay $\mathbf{S}$-module, then  $i'_\lambda N$ is a Cohen-Macaulay $\mathbf{S}'$-module;
\item the pair $(i'_\lambda, i')$ is adjoint.
\end{enumerate}
\end{lem}

\begin{proof}
The statements (1)-(5) are proved with a similar argument as in Lemma \ref{lem:i'}.
We point out that the statement (3) is essentially contained in \cite[Proposition 9.4]{GL90}, while
the  statement (6) is implicitly stated in \cite[p.325]{GL90}.

Here we indicate  the construction of the isomorphism $\Phi \colon {\rm Hom}_{{\rm mod}^{\mathbf{L}'}\; \mathbf{S}'}(i'_\lambda N, M)\simeq {\rm Hom}_{{\rm mod}^{\mathbf{L}}\; \mathbf{S}}(N, i'M)$ for
the adjoint pair in (6). It
sends $f\colon i'_\lambda N \rightarrow M$ to $\Phi(f)$ such that its restriction on
$N_{\vec{l}}$, that is, $N_{\vec{l}}\rightarrow (i'M)_{\vec{l}}$, is restricted from $f$, except
that in the case $l_n< p_n-p'_n$, it is given by $N_{\vec{l}}\stackrel{x_n^{p_n-p'_n-l_n}} \rightarrow
N_{\vec{l}+(p_n-p'_n-l_n)\vec{x}_n}=(i'_\lambda N)_{\phi'^{-1}(\vec{l}-l_n\vec{x}_n)} \rightarrow M_{\phi'^{-1}(\vec{l}-l_n\vec{x}_n)}=(i'M)_{\vec{l}}$,
where the second map is restricted from $f$.
\end{proof}

We define a functor  $i'_\rho \colon {\rm mod}^{\mathbf{L}}\; \mathbf{S}\rightarrow {\rm mod}^{\mathbf{L}'}\; \mathbf{S}'$ as follows.
For an $\mathbf{L}$-graded $\mathbf{S}$-module $N$, we set $i'_\rho N=\oplus_{\vec{l}\in \mathbf{L}'} (i'_\rho N)_{\vec{l}}$ such that  $(i'_\rho N)_{\vec{l}}=N_{\phi'(\vec{l})}$
if $l_n=0$, and $(i'_\rho N)_{\vec{l}}=N_{\phi'(\vec{l})+(p_n-p'_n)\vec{x}_n}$, otherwise. The action of $u$, $v$ and $x_i$ on $i'_\lambda N$ is induced by the one on $N$, except that $x_n$ acts by
$x_n^{p_n-p'_n}$ (in $\mathbf{S}$) on $(i'_\rho N)_{\vec{l}}$  provided that $l_n=0$. The obtained $\mathbf{L}'$-graded $\mathbf{S}'$-module $i'_\rho N$ is finitely generated by Lemma \ref{lem:f.g.}.
The action of $i'_\rho$ on morphisms is defined naturally. We observe that $(i'_\rho N)(\vec{x}_n)=i'_\lambda(N(\vec{x}_n))$
for any $\mathbf{L}$-graded $\mathbf{S}$-module $N$.

The following lemma is dual to Lemma \ref{lem:leftadj}.

\begin{lem}\label{lem:rightadj}
Use the notation above. Then the  following statements hold:
\begin{enumerate}
\item the functor $i'_\rho$ is exact;
\item $i'_\rho (N(\vec{x}_i))=(i'_\rho N)(\vec{x}_i)$ for any graded $\mathbf{S}$-module $N$ and $1\leq i< n$;
\item $i'_\rho (\mathbf{S}(\vec{l}))\simeq \mathbf{S}'(\phi'^{-1}(\vec{l}))$  if $0\leq l_n< p'_n$, and $i'_\rho (\mathbf{S}(\vec{l}))\simeq \mathbf{S}'(\phi'^{-1}(\vec{l}-(l_n-p'_n+1) \vec{x}_n))$, otherwise;
\item $i'_\rho (k(\vec{l}))=k(\phi'^{-1}(\vec{l}))$ if $0\leq l_n< p'_n$, and $i'_\rho (k(\vec{l}))=0$, otherwise;
\item if $N$ is a Cohen-Macaulay $\mathbf{S}$-module, then  $i'_\rho N$ is a Cohen-Macaulay $\mathbf{S}'$-module;
\item the pair $(i', i'_\rho)$ is adjoint. \hfill $\square$
\end{enumerate}
\end{lem}

We have built adjoint pairs $(i'_\lambda, i')$ and $(i', i'_\rho)$ of exact functor on graded modules.
Observe that these functors preserve finite dimensionality. By abuse of notation we have the induced
functor $i'\colon {\rm coh}\; \mathbb{X}'\rightarrow {\rm coh}\; \mathbb{X}$ which has a left
adjoint $i'_\lambda\colon {\rm coh}\; \mathbb{X}\rightarrow {\rm coh}\; \mathbb{X}'$ and a right adjoint
$i'_\rho \colon {\rm coh}\; \mathbb{X}\rightarrow {\rm coh}\; \mathbb{X}'$; see Lemma \ref{lem:quotientAbel}.
These induced functors are all exact, and $i'$ is fully faithful.

Recall that the sheafification functor ${\rm mod}^\mathbf{L}\; \mathbf{S}\rightarrow {\rm coh}\; \mathbb{X}$
induces an equivalence ${\rm CM}^\mathbb{L}(\mathbf{S})\simeq {\rm vect}\; \mathbb{X}$ of exact categories, which
identifies projective modules with line bundles; see Lemma \ref{lem:MCM-vb}. Here, the exact structure on ${\rm vect}\; \mathbb{X}$ is given by the distinguished exact sequences. It follows then that the obtained three functors on sheaves restricts to three
exact functors on the category of vector bundles. Moreover, these restricted functors preserve line bundles; see Lemmas \ref{lem:i'}(3), \ref{lem:leftadj}(3) and \ref{lem:rightadj}(3). Therefore, these restricted functors induces triangle
functors on the stable category of vector bundles. Applying Lemma \ref{lem:stable} we obtain two  adjoint pairs
$(i'_\lambda, i')$ and $(i', i'_\rho)$ of triangle functors, where the triangle functor $i'\colon \underline{\rm vect}\; \mathbb{X}' \longrightarrow \underline{\rm vect}\; \mathbb{X}$ is fully faithful. Here we abuse the notation again.

We have the following immediate consequence of Lemma \ref{lem:recollement1}; compare \cite[Theorem 4.3.1]{CK}.

\begin{prop}\label{prop:recollement}
Keep the notation as above. Then we have the following recollement
\[\xymatrixrowsep{3pc} \xymatrixcolsep{2pc}\xymatrix{
  \underline{\rm vect}\; \mathbb{X}'  \;\ar[rr]|-{i'} &&\;  \underline{\rm vect}\; \mathbb{X} \; \ar[rr]|-{q}
\ar@<1.2ex>[ll]^-{i'_\rho} \ar@<-1.2ex>[ll]_-{i'_\lambda} &&
\; \underline{\rm vect}\; \mathbb{X}/{{\rm Im}\; i'}, \ar@<1.2ex>[ll]^-{}\ar@<-1.2ex>[ll]_-{}
}\]
where $q\colon \underline{\rm vect}\; \mathbb{X} \rightarrow \underline{\rm vect}\; \mathbb{X}/{{\rm Im}\; i'}$ denotes the quotient
functor. \hfill $\square$
\end{prop}

In general, we do not know much about the Verdier quotient category  $\underline{\rm vect}\; \mathbb{X}/{{\rm Im}\; i'}$
in the recollement above. Note that the case $n=2$ is boring, since then the three triangulated categories in the recollement
 are trivial. We will see that if $n=3$, that is, the weight sequence of the weighted projective line $\mathbb{X}$
 has length $3$, then we have an explicit description of the quotient category.

\section{The main result}

In this section we describe the quotient category appearing in the recollement of Proposition \ref{prop:recollement}
under the condition that the weight sequence has length $3$. This yields our main result, where an explicit
recollement consisting of the stable categories of vector bundles is given; see Theorem \ref{thm:maintheorem}.

Let $\mathbf{p}=(p_1, p_2, p_3)$ be a weight sequence of length $3$, and let $\boldsymbol \lambda=(\lambda_1,\lambda_2, \lambda_3)$
be a parameter sequence. Denote by $\mathbb{X}=\mathbb{X}(\mathbf{p}, \boldsymbol \lambda)$ the corresponding weighted
projective line. Note that  the category of
coherent sheaves on the  weighted projective line $\mathbb{X}$, up to equivalence,  does not depends on the choice of the parameter sequence $\boldsymbol \lambda$, since the weight sequence has length $3$; compare \cite[Proposition 9.1]{GL90}.
For this reason, as we do in Introduction, the weighted projective line $\mathbb{X}(\mathbf{p}, \boldsymbol \lambda)$ is sometimes written as $\mathbb{X}(\mathbf{p})$.

 Fix a positive integer $p'_3$ such that $p'_3\leq p_3$. Set $p''_3=p_3-p'_3+1$. Set $\mathbf{p}'=(p_1, p_2, p'_3)$ and  $\mathbb{X}'=\mathbb{X}(\mathbf{p}', \boldsymbol \lambda)$. Recall that $\mathbf{S}'=\mathbf{S}(\mathbf{p}', \boldsymbol \lambda)$, is the homogeneous coordinate ring of $\mathbb{X}'$, which is graded by $\mathbf{L}'=\mathbf{L}(\mathbf{p}')$.  Similarly we have the notation $\mathbf{p}''$, $\mathbb{X}''$, $\mathbf{S}''$ and $\mathbf{L}''$.

Recall from Section 4 the explicitly given exact functor $i'\colon {\rm mod}^{\mathbf{L}'}\; \mathbf{S}'\rightarrow
{\rm mod}^\mathbf{L}\; \mathbf{S}$, which allows an exact left adjoint $i'_\lambda$ and an exact right adjoint $i'_\rho$.
Observe that all these functors preserve projective modules. These exact functors extend to triangle functors between the corresponding  bounded derived categories, which preserve perfect complexes. Applying a triangulated analogue of Lemma \ref{lem:quotientAbel}, we get the induced triangle functor $i'\colon \mathbf{D}_{\rm sg}^{\mathbf{L}'}(\mathbf{S}')\rightarrow  \mathbf{D}_{\rm sg}^{\mathbf{L}}(\mathbf{S})$, which allows a left adjoint $i_\lambda$ and a right adjoint $i'_\rho$.
We have observed in Section 4 that the three exact functors on module categories induces the corresponding
triangle functors between the stable categories of vector bundles.

We recall the triangle equivalence in Proposition \ref{prop:vb-singularity}, which will be
denoted by $F\colon \underline{\rm vect}\; \mathbb{X}\rightarrow
\mathbf{D}_{\rm sg}^\mathbf{L}(\mathbf{S})$. Similarly, we have a triangle equivalence $F'\colon \underline{\rm vect}\; \mathbb{X}'\rightarrow
\mathbf{D}_{\rm sg}^{\mathbf{L}'}(\mathbf{S}')$.

 The following immediate observation states that these triangle equivalences are compatible with
 the functors $i'$, $i'_\lambda$ and $i'_\rho$ defined on both sides.

\begin{lem}\label{lem:compatible}
Keep the notation as above. Then we have natural isomorphisms $i'F'\simeq Fi'$,  $i'_\lambda F\simeq F' i'_\lambda $ and $i'_\rho F\simeq F' i'_\rho$. \hfill $\square$
\end{lem}

Recall that $\mathbf{S}''=\mathbf{S}(\mathbf{p}'', \boldsymbol \lambda)$ and $\mathbb{X}''=\mathbb{X}(\mathbf{p}'', \boldsymbol \lambda)$. Then we have the exact fully faithful functor $i''\colon {\rm mod}^{\mathbf{L}''}\; \mathbf{S}''\rightarrow {\rm mod}^{\mathbf{L}}\; \mathbf{S}$ which admits an exact left adjoint $i''_\lambda$ and an exact right adjoint $i''_\rho$; see Section 4. These functors induce the corresponding functors on  the stable
categories of vector bundles and the graded singularity categories; these induced functors
are still denoted by $i''$, $i''_\lambda$ and $i''_\rho$.  Moreover, the triangle equivalences $F$ and  $F''\colon \underline{\rm vect}\; \mathbb{X}''\rightarrow \mathbf{D}_{\rm sg}^{\mathbf{L}''}(\mathbf{S}'')$ are compatible with these functors; compare Lemma \ref{lem:compatible}.

We are in the position to state and prove our main result.

\begin{thm}\label{thm:maintheorem}
Keep the assumption and notation as above. Then we have the following recollement of triangulated categories
\[\xymatrixrowsep{3pc} \xymatrixcolsep{2pc}\xymatrix{
  \underline{\rm vect}\; \mathbb{X}'  \;\ar[rrr]|-{i'} &&& \;  \underline{\rm vect}\; \mathbb{X} \; \ar[rrr]|-{i''_\lambda ((1-p'_3)\vec{x}_3)}
\ar@<1.3ex>[lll]^-{i'_\rho} \ar@<-1.3ex>[lll]_-{i'_\lambda} &&&
\; \underline{\rm vect}\; \mathbb{X}'' . \ar@<1.3ex>[lll]^-{((p'_3-1)\vec{x}_3)i''} \ar@<-1.3ex>[lll]_-{(p'_3\vec{x}_3)i'' (-\vec{x}_3)}
}\]
In particular, we have a triangle equivalence $\underline{\rm vect}\; \mathbb{X}/{{\rm Im}\; i'}\simeq \underline{\rm vect}\; \mathbb{X}''$.
\end{thm}

\begin{proof}
Set $j''=i'_\lambda ((1-p'_3)\vec{x}_3)$, $j''_\lambda=(p'_3\vec{x}_3)i'' (-\vec{x}_3)$ and $j''_\rho=((p'_3-1)\vec{x}_3)i''$.
Recall that $(i'_\lambda, i')$ and $(i', i''_\rho)$ are adjoint pairs. Similarly, $(i''_\lambda, i'')$ and $(i'', i''_\rho)$ are
adjoint pairs. Then it follows that $(j'', j''_\rho)$ is an adjoint pair. Note that $j''=(\vec{x}_3) i''_\rho (-p'_3\vec{x}_3)$, since
we have $(\vec{x}_3)i''_\rho=i''_\lambda (\vec{x}_3)$; see Section 4. Then we have that  $(j''_\lambda, j'')$ is also an adjoint pair. Then the above diagram satisfies the conditions (R1) and (R2).  We will apply Lemma \ref{lem:recollement2}. Then it suffices to show
that $j''i'\simeq 0$ and ${\rm thick}\langle {\rm Im}\; i'\cup {\rm Im}\; j''_\lambda\rangle =\underline{\rm vect}\; \mathbb{X}$.

Recall that the triangle equivalences $F$, $F'$ and $F''$ are compatible with the degree-shift functors and the twist functors, and also with  the six functors $i'$, $i'_\lambda$, $i'_\rho$, $i''$, $i''_\lambda$ and $i''_\rho$; see Proposition \ref{prop:vb-singularity} and Lemma \ref{lem:compatible}. Then it follows that they are compatible with the functors $j''$, $j''_\lambda$ and $j''_\rho$. Using these three triangle equivalences again, it suffices to show the following two statements: (1) the composite $\mathbf{D}_{\rm sg}^{\mathbf{L}'}(\mathbf{S}')\stackrel{i'}\rightarrow \mathbf{D}_{\rm sg}^\mathbf{L}(\mathbf{S})\stackrel{j''}\rightarrow \mathbf{D}_{\rm sg}^{\mathbf{L}''}(\mathbf{S}'')$ is zero; (2) the union of the  images of the functors $i'\colon \mathbf{D}_{\rm sg}^{\mathbf{L}'}(\mathbf{S}')\rightarrow \mathbf{D}_{\rm sg}^\mathbf{L}(\mathbf{S})$ and $j''_\lambda \colon \mathbf{D}_{\rm sg}^{\mathbf{L}''}(\mathbf{S}'')\rightarrow \mathbf{D}_{\rm sg}^\mathbf{L}(\mathbf{S})$ generates $\mathbf{D}_{\rm sg}^\mathbf{L}(\mathbf{S})$.

Recall from Lemma \ref{lem:singularity}(2) that the category $\mathbf{D}_{\rm sg}^{\mathbf{L}'}(\mathbf{S}')$  is generated by
$\{qk(\vec{l})\; |\; \vec{l}\in \mathbf{L}'\}$, where $q$ is the quotient functor. To see the statement (1) it suffices to show that $j''i'(qk(\vec{l}))\simeq 0$ for each $\vec{l}$ in $\mathbf{L}'$. We write $\vec{l}$ in this normal form. Then we observe that $j''i'(k(\vec{l}))=0$ (as an $\mathbf{S}''$-module) if $l_3>0$. If $l_3=0$, we have $i'(k(\vec{l}))=\mathbf{S}(\phi'(\vec{l}))/(x_1, x_2, x_3^{p_3-p'_3+1})=\mathbf{S}(\phi'(\vec{l}))/(x_1, x_2, x_3^{p''_3})$.  Here, we recall
that $n=3$, that is, the weight sequence of $\mathbb{X}$ has length $3$. Then we have $j''i'(k(\vec{l}))=\mathbf{S}''(\phi''^{-1}(\phi'(\vec{l}))+\vec{c})/{(x_1, x_2, x_3^{p''_3})}=\mathbf{S}''(\phi''^{-1}(\phi'(\vec{l}))+\vec{c})/(x_1, x_2)$. Since $\{x_1, x_2\}$ is a
(homogeneous) regular sequence in $\mathbf{S}''$, the $\mathbf{S}''$-module $j''i'(k(\vec{l}))$ has finite projective dimension.
Hence $j''i'(qk(\vec{l}))=q(j''i'k(\vec{l}))\simeq 0$ in $\mathbf{D}_{\rm sg}^{\mathbf{L}''}(\mathbf{S}'')$.

 It remains to show the statement (2). We write an element  $\vec{l}\in \mathbf{L}$ in its normal form. We observe that by Lemma \ref{lem:i'}(4) that $qk(\vec{l})$ lies in the image of $i'\colon \mathbf{D}_{\rm sg}^{\mathbf{L}'}(\mathbf{S}')\rightarrow \mathbf{D}_{\rm sg}^\mathbf{L}(\mathbf{S})$
 provided that $0<l_3<p'_3$. Similarly $qk(\vec{l})$ lies in the image of $j''_\lambda \colon \mathbf{D}_{\rm sg}^{\mathbf{L}''}(\mathbf{S}'')\rightarrow \mathbf{D}_{\rm sg}^\mathbf{L}(\mathbf{S})$ provided that $p'_3+1\leq l_3< p_3$.
 Hence we have that $qk(\vec{l})$ lies in ${\rm thick}\langle {\rm Im}\; i'\cup {\rm Im}\; j''_\lambda\rangle$ provided that $l_3>0$.
 We will show that  $qk(\vec{l})$ lies in ${\rm thick}\langle {\rm Im}\; i'\cup {\rm Im}\; j''_\lambda\rangle$ in the case
 $l_3=0$. Then by Lemma \ref{lem:singularity}(2) we are done.

  We assume that $l_3=0$. By Lemma \ref{lem:i'}(4) we have a short
 exact sequence in ${\rm mod}^\mathbf{L}\; \mathbf{S}$
 $$0\longrightarrow K \longrightarrow i'(k(\phi'^{-1}(\vec{l})))\longrightarrow k(\vec{l})\longrightarrow 0,$$
 where $K$ is a finite dimension module with composition factors $\{k(\vec{l}-\vec{x}_3), k(\vec{l}-2\vec{x}_3), \cdots, k(\vec{l}-p''_3\vec{x}_3)\}$. This exact sequence induces a triangle in $\mathbf{D}_{\rm sg}^{\mathbf{L}}(\mathbf{S})$.
 Observe that $qK$ lies in ${\rm thick}\langle k(\vec{l}-\vec{x}_3), k(\vec{l}-2\vec{x}_3), \cdots, k(\vec{l}-p''_3\vec{x}_3) \rangle $, and thus by above in ${\rm thick}\langle {\rm Im}\; i'\cup {\rm Im}\; j''_\lambda\rangle$. Then the induced triangle forces
 that $qk(\vec{l})$ lies in  ${\rm thick}\langle {\rm Im}\; i'\cup {\rm Im}\; j''_\lambda\rangle$, completing the proof.
\end{proof}

\begin{rem}
The above proof yields the following two recollements, both of which are isomorphic to the recollement
above. Here we use the equivalences in Lemma \ref{lem:MCM-vb} and Proposition \ref{prop:vb-singularity}.
\[\xymatrixrowsep{3pc} \xymatrixcolsep{2pc}\xymatrix{
 \mathbf{D}_{\rm sg}^{\mathbf{L}'}(\mathbf{S}')  \;\ar[rrr]|-{i'} &&& \;  \mathbf{D}_{\rm sg}^{\mathbf{L}}(\mathbf{S}) \; \ar[rrr]|-{i''_\lambda ((1-p'_3)\vec{x}_3)}
\ar@<1.3ex>[lll]^-{i'_\rho} \ar@<-1.3ex>[lll]_-{i'_\lambda} &&&
\; \mathbf{D}_{\rm sg}^{\mathbf{L}''} (\mathbf{S}''). \ar@<1.3ex>[lll]^-{((p'_3-1)\vec{x}_3)i''} \ar@<-1.3ex>[lll]_-{(p'_3\vec{x}_3)i'' (-\vec{x}_3)}
}\]

\[\xymatrixrowsep{3pc} \xymatrixcolsep{2pc}\xymatrix{
  \underline{\rm CM}^{\mathbf{L}'}(\mathbf{S}')  \;\ar[rrr]|-{i'} &&& \;   \underline{\rm CM}^{\mathbf{L}}(\mathbf{S}) \; \ar[rrr]|-{i''_\lambda ((1-p'_3)\vec{x}_3)}
\ar@<1.3ex>[lll]^-{i'_\rho} \ar@<-1.3ex>[lll]_-{i'_\lambda} &&&
\; \ \underline{\rm CM}^{\mathbf{L}''}(\mathbf{S}''). \ar@<1.3ex>[lll]^-{((p'_3-1)\vec{x}_3)i''} \ar@<-1.3ex>[lll]_-{(p'_3\vec{x}_3)i'' (-\vec{x}_3)}
}\]

\end{rem}

\vskip 10pt

\noindent {\bf Acknowledgements} \;  This research was mainly carried out during the author's visit at the University of Bielefeld with a support by Alexander von Humboldt Stiftung. He would like to thank Professor  Henning Krause and the faculty of Fakult\"{a}t f\"{u}r Mathematik for their hospitality.

\bibliography{}

\begin{thebibliography}{999}


%
%
\bibitem{BBD} {\sc A A. Beilinson, J. Bernstein and P. Deligne},
Faisceaux Perves, Ast\'{e}rique {\bf 100}, Soc. Math. France, 1982.

%
%
%
%
%
%
%
%
%


%

\bibitem{BK}{\sc A.I. Bondal and M.M. Kapranov,} {\em Representable functors, Serre
functors, and reconstructions,} Izv. Akad. Nauk SSSR Ser. Mat. {\bf
53} (1989), 1183--1205.

%
%


\bibitem{Buc87} {\sc R.O. Buchweitz,} Maximal Cohen-Macaulay Modules
and Tate Cohomology over Gorenstein Rings. Unpublished manuscript,
1987.

%
%
%
%
%
%
%
%
%
%
%




\bibitem{Ch11} {\sc X.W. Chen}, {\em Three results on Frobenius
categories}, Math. Z., accepted; arXiv:1004.4540v1.


\bibitem{Ch11'} {\sc X.W. Chen}, {\em Unifying two results of Orlov on
singularity categories}, Abh. Math. Semin. Univ. Hambg. {\bf 80} (2010), 207--212.



\bibitem{CK} {\sc X.W. Chen and H. Krause}, {\em Expansions of abelian categories},
arXiv:1009.3456v1.



\bibitem{CL} {\sc X.W. Chen and S. Ladkani}, {\em The $F$-inflation category and its stable
category}, in preparation.



%
%
%
%
%
%
%



\bibitem{Ga62}{\sc P. Gabriel}, {\em Des cat\'egories
ab\'eliennes,} Bull. Soc. Math. France {\bf 90} (1962), 323--448.



\bibitem{GL87} {\sc W. Geigel and H. Lenzing,} {\em A class of weighted projective curves arising in
representation theory of finite dimensional algebras,} in:
Singularities, representations of algebras and vector bundles,
Lecture Notes in Math. {\bf 1273}, 265--297, Springer, 1987.

%
%

\bibitem{GL90} {\sc W. Geigle and H. Lenzing,} {\em Perpendicular categories with applications
to representations and sheaves,} J. Algebra {\bf 144} (1991),
273--343.

%
%

\bibitem{Ha88}{\sc D. Happel,} Triangulated Categories in the
Representation Theory of Finite Dimensional Algebras.  London Math.
Soc. Lecture Notes Ser. {\bf 119}, Cambridge Univ. Press, Cambridge,
1988.

%
%
%
%
%
%

\bibitem{Ke90} {\sc B. Keller,} {\em Chain complexes and stable
categories,} Manuscripta Math. {\bf 67} (1990), 379-417.

%
%

\bibitem{Ke96} {\sc B. Keller,} {\em Derived categories and their
uses,} Handbook of Algebra {\bf 1}, 671--701, North-Holland,
Amsterdam, 1996.

%
%
%
%
%

\bibitem{KMV08} {\sc B. Keller, D. Murfet and M. Van den Bergh, } {\em On two examples by Iyama and
Yoshino}, Compositio Math., to appear, arXiv: 0803.0720v3.



\bibitem{KLM1} {\sc D. Kussin, H. Lenzing and H. Meltzer}, {\em Triangle singularities, ADE-chains
and weighted projective lines}, preprint, 2008.

%
\bibitem{KLM2} {\sc D. Kussin, H. Lenzing and H. Meltzer}, {\em Nilpotent operators and
weighted projective lines}, arXiv:1002.3797v1.
%
%

\bibitem{KLM3} {\sc D. Kussin, H. Lenzing and H. Meltzer}, {\em Weighted projective lines and
invariant flags of nilpotent operators}, in preparation.

%
%
%
%

\bibitem{Lende} {\sc H. Lenzing and J. de la Pena}, {\em Spectral analysis of finite dimensional algebras and singularities},
In: Trends in representation theory of algebras and related topics, 541¨C-588, EMS Ser. Congr. Rep., Eur. Math. Soc.,
 Z\"{u}rich, 2008.


%
%
%
%

\bibitem{MacL} {\sc S. Mac Lane}, Categories for the Working
Mathematicians, Grad. Text in Math. {\bf 5}, Springer, 1998.

%
%
%
%
%

\bibitem{Or04} {\sc D. Orlov}, {\em Triangulated categories of singularities and D-branes in Landau-Ginzburg
models}, Trudy Steklov Math. Institute {\bf 204} (2004), 240--262.


\bibitem{Or09}{\sc D. Orlov}, {\em Derived categories of coherent sheaves and triangulated categories
of singularities,} in: {\it Algebra, arithmetic, and geometry: in honor of Yu. I. Manin. Vol. II},
503--531, Progr. Math. {\bf 270}, Birkh\"{a}user Boston, Inc., Boston, MA, 2009.


\bibitem{Or10} {\sc D. Orlov}, {\em Formal completions and idempotent completions of triangulated
categories of singularities}, Adv. Math., to appear, arXiv: 0901.1859v1.



%
%
%
%

\bibitem{Qui73} {\sc D. Quillen,} {\em Higher algebraical K-theory
I}, Springer Lecture Notes in Math. {\bf 341}, 1973, 85--147.

%
%
%
%

\bibitem{RS06} {\sc C.M. Ringel and M. Schmidmeier}, {\em Submodule categories of wild representation type},
J. Pure and Applied Algebra {\bf 205} (2) (2006), 412--422.

%
%
%
%

\bibitem{RS08}{\sc C.M. Ringel and M. Schmidmeier,} {\em The Auslander-Reiten translation in
submodule category}, Trans. Amer. Math. Soc. {\bf 360} (2) (2008),
691--716.

%
%

\bibitem{RS08'} {\sc C.M. Ringel and M. Schmidmeier,} {\em Invariant subspaces of
nilpotent linear operators}, J. Reine Angew. Math. {\bf 614} (2008),
1--52.

%
%
%
%


\bibitem{Ver77} {\sc J.L. Verdier}, {\em Categories derive\'{e}es},
in SGA 4 1/2, Lecture Notes in Math. {\bf 569}, Springer, Berlin,
1977.







\end{thebibliography}

\vskip 10pt

{\footnotesize \noindent Xiao-Wu Chen, Department of Mathematics,
University of Science and Technology of
China, Hefei 230026, P. R. China \\
Homepage: http://mail.ustc.edu.cn/$^\sim$xwchen}

\end{document}